\documentclass{amsart}
\usepackage{enumerate}
\usepackage{colortbl}
\newtheorem{theorem}{Theorem}[section]
\newtheorem{lemma}[theorem]{Lemma}
\newtheorem{corollary}[theorem]{Corollary}
\newtheorem{proposition}[theorem]{Proposition}
\theoremstyle{definition}
\newtheorem{definition}[theorem]{Definition}
\newtheorem{example}[theorem]{Example}

\theoremstyle{remark}
\newtheorem{remark}[theorem]{Remark}

\numberwithin{equation}{section}

\def\C{\mathbb C}
\def\R{\mathbb R}
\def\X{\mathbb X}
\def\Z{\mathbb Z}
\def\Y{\mathbb Y}
\def\Z{\mathbb Z}
\def\N{\mathbb N}

\begin{document}

\title[Polynomially bounded sequences]{
A Spectral Theory of Polynomially Bounded Sequences and Applications to the Asymptotic Behavior of Discrete Systems}

\author[N.V. Minh]{Nguyen Van Minh}
\address{Department of Mathematics and Statistics, University of Arkansas at Little Rock, 2801 S University Ave, Little Rock, AR 72204. USA}
\email{mvnguyen1@ualr.edu}

\author[H. Matsunaga]{Hideaki Matsunaga}
\address{Department of Mathematical Sciences, Osaka Prefecture University,
Sakai 599-8531, Japan}
\email{hideaki@ms.osakafu-u.ac.jp} 

\author[N.D. Huy]{Nguyen Duc Huy}
\address{VNU-University of Education, Vietnam National University, Hanoi; 144 Xuan Thuy, Cau Giay, Hanoi, Vietnam}
\email{huynd@vnu.edu.vn}

\author[V.T Luong]{Vu Trong Luong}
\address{VNU-University of Education, Vietnam National University, Hanoi; 144 Xuan Thuy, Cau Giay, Hanoi, Vietnam}
\email{vutrongluong@gmail.com}

\thanks{This work is supported by Viet Nam Ministry of Education and Training under grant number B2019-TTB-01 and 
JSPS KAKENHI Grant Number JP19K03524.}
\thanks{The authors thank the anonymous referee for his comments on the previous version of the paper for us to improve the presentation of the paper as well as to correct several errors.}

\date{\today}

\begin{abstract} 
In this paper using a transform defined by the translation operator we introduce the concept of spectrum of sequences  that are bounded by $n^\nu$, where $\nu$ is a natural number. 
We apply this spectral theory to study the asymptotic behavior of solutions of fractional difference equations of the form $\Delta^\alpha x(n)=Tx(n)+y(n)$, $n\in \N$, where $0<\alpha\le 1$. One of the obtained results is an extension of a famous Katznelson-Tzafriri Theorem, saying that if the $\alpha$-resolvent operator $S_\alpha$ satisfies $\sup_{n\in\N} \| S_\alpha (n)\| /n^\nu <\infty$ and for all $z_0\in \{z\in \C: \ |z|=1\}$, but $z_0=1$, the complex function 
$(z^{1-\alpha}(z-1)^\alpha -T)^{-1}$ \ exists and is holomorphic in a neighborhood of $z_0$,
then
\begin{align*}
\lim_{n\to \infty} \frac{1}{n^\nu} \sum_{k=0}^{\nu+1} \frac{(\nu+1)!}{k!(\nu+1-k)!} (-1)^{\nu+1+k} S_\alpha (n+k) =0.
\end{align*}
Three concrete examples are also included to illustrate the obtained results. 
\end{abstract}

\maketitle

\section{Introduction}
Let us consider difference equations of the form
\begin{align}\label{1.1}
   x(n+1)=Tx(n)+y(n), \quad n\in \N,
\end{align}
where $T$ is a bounded operator in a Banach space $\X$, $\{x(n)\}_{n=1}^\infty$ and $\{ y(n)\}_{n=1}^\infty$ are sequences in $\X$. The asymptotic behavior of solutions of the above mentioned equations is a central topic in Analysis and Dynamical Systems. There are numerous methods for this study of this topic.  The reader is referred to \cite{ela} and its references for information on classical methods of Dynamical Systems in the finite dimensional case. On the other hand, in the infinite dimensional case, by Harmonic Analysis and Operator Theory, many results on the asymptotic behavior of solutions of Eq.~(\ref{1.1}) have been obtained, see e.g. \cite{aba, arebathieneu, baspry, bas, batyea, debvin, kattza, min, sei, vu1, vu2}. Among many interesting results in this direction is a famous theorem due to Katznelson-Tzafriri (see \cite{kattza}) saying that if $T$ is a bounded operator in a Banach space $\X$ such that
\begin{align}\label{1.2}
   \sup_{n\in\N}\| T^n\| <\infty ,
\end{align}
and $\sigma (T)\subset \{1\}$, then
\begin{align}
   \lim_{n\to \infty}(T-I)T^n=0.
\end{align}
There are a lot of extensions and improvements of this result as well as simple proofs of it, see e.g. \cite{aba, arebathieneu, bas, debvin, kattza, min, sei, vu1, vu2} and the references therein. 

As shown in \cite{vu1} the above mentioned Katznelson-Tzafriri Theorem is equivalent to its weaker version for individual orbits. Namely, the following statement: Let $T$ be a bounded operator in a Banach space $\X$ such that (\ref{1.2}) holds and $\sigma (T)\subset \{1\}$. Then, for each $x\in \X$
\begin{align}
   \lim_{n\to\infty} (T-I)T^nx=0.
\end{align}
In \cite{min} a simple proof of this weaker version is given, based on a transform associated with the translation operator of sequences.

The main concern of this paper is to extend the above mentioned Katznelson-Tzafriri Theorem to fractional difference equations of the form
\begin{align}
   \Delta^\alpha x(n)=Tx(n)+y(n), \quad n\in \N,
\end{align}
where $0<\alpha \le 1$, the operator $\Delta^\alpha$ (the fractional difference operator in the sense of Riemann-Liouville) and other operators are defined as follows (see \cite{atielo,liz} and its references for more details): for each $n\in \N$,
\begin{align}
   (\Delta^\alpha )f(n) &=  \Delta ^1 \circ      \nabla_0^{-(1-\alpha)}f(n), \nonumber \\
   (\Delta^1 f)(n) &= f(n+1)-f(n),  \nonumber\\
   (\nabla_0^{-\alpha} )f(n)&= \sum_{k=0}^n k^\alpha (n-k)f(k), \nonumber \\
   k^\alpha (j)&= \frac{\Gamma (\alpha +j)}{\Gamma (\alpha )\Gamma (j+1)}, \label{k alpha}
\end{align}
where $\Gamma (\cdot )$ is the gamma function  defined below.
Our method relies on a spectral theory of polynomially bounded sequences that will be presented in the next sections, and that would be of independent interest. The obtained results will be illustrated in simple cases of ordinary difference equations and then stated for fractional difference equations. Our main result is Theorem \ref{the kat for frac}. To our best knowledge, it is a new extension of the Katznelson-Tzafriri Theorem to fractional difference equations.

\section{Preliminaries and Notations}
\subsection{Notations}
Throughout this paper we will denote by $\N$, $\Z$, $\R$, $\C$ the set of natural numbers, integers, real numbers and the complex plane, respectively. For $z\in \C$, $\Re \hspace{0.3mm} z$ stands for its real part. 
The gamma function $\Gamma (z)$ is defined to be
\begin{align*}
   \Gamma (z)=\int^\infty _0 x^{z-1}e^{-x} dx , \quad \Re \hspace{0.3mm} z >0.
\end{align*}
In this paper $\X$ is assumed to be a complex Banach space. For an element $v\in\X$, $\| v\|$ denotes the norm of $v$. $L(\X)$ denotes the space of all bounded linear operators from $\X$ to itself. We also use  the following standard notations: $\rho (T)$ denotes the resolvent set of a given operator $T$, that is, $\rho(T):=\{ \lambda \in \C :\ (\lambda -T)^{-1} \mbox{exists}\}$, and $\sigma (T):= \C \backslash \rho (T)$. For each $\lambda \in \rho (T)$ we denote $R(\lambda , T):=(\lambda -T)^{-1}$. 
Moreover, we will denote by $\Gamma$ the unit circle in the complex plane. 

For a given nonnegative integer $\nu$, we denote by $l_\infty^\nu (\X)$ the space of all sequences in a Banach space $\X$ such that
\begin{align*}
   \sup_{n\in \N} \frac{\| x(n)\|}{n^\nu } < \infty .
\end{align*}
It is easy to see that $l_\infty^\nu (\X)$ is a Banach space with norm
\begin{align*}
   \| x\|_\nu  := \sup_{n\in \N}  \frac{\| x(n)\|}{n^\nu }
\end{align*}
for each $x=\{ x(n)\}_{n\in \N}$.
We will denote by $c_0^\nu (\X)$ the subspace of  $l_\infty^\nu (\X)$ consisting of all sequences $\{x(n)\}_{n\in\N}$ such that
\begin{align*}
   \lim_{n\to \infty} \frac{\| x(n)\|}{n^\nu }=0.
\end{align*}

We can check that $c_0^\nu (\X)$ is a complete subspace of $l_\infty^\nu (\X)$, so the quotient space
\begin{align*}
   \Y:= l_\infty^\nu (\X)/c_0^\nu (\X)
\end{align*}
is well defined as a Banach space. If $f\in l_\infty^\nu (\X)$ we will denote its equivalence class by $\bar f$.
In the space $l_\infty^\nu (\X)$ let us consider the translation operator $S$ defined as
\begin{align*}
   [Sx](n)=x(n+1), \quad n\in \N, \enskip x\in l_\infty^\nu  (\X).
\end{align*}
This is a bounded operator. Moreover, this operator leaves $c_0^\nu (\X)$ invariant. Hence, it induces an operator $\bar S$ on $\Y$.

\subsection{Vector-valued holomorphic functions} 
In this paper we say that a function $f(z)$ defined for all $z\in \Omega \subset \C$ with values in a complex Banach space $\X$ is holomorphic (or analytic) for  $z\in \Omega$ if for each $z_0\in \Omega$
\begin{align*}
   f'(z_0):=\lim_{h\to 0, h\not= 0} \frac{f(z_0+h)-f(z_0)}{h}
\end{align*}
exists. 
A family of continuous functionals $W \subset \X^*$ is said to be {\it separating} if $x\in \X$ and $\langle x,\phi\rangle=0$ for all $\phi \in W$, then $x=0$. We will need the following whose proof can be found in  \cite[Theorem 3.1]{arenik}, or \cite[Theorem A.7]{arebathieneu}:
\begin{theorem}\label{the anal}
Let $\Omega \subset \C$ be open and connected, and let $f:\Omega \to \X$ be bounded on every compact subset of $\Omega$. Assume further that $W\subset \X^*$ is separating subset such that $x^* \circ f$ is holomorphic for all $x^*\in W$. Then $f$ is holomorphic.
\end{theorem}

We will need an auxiliary result that is a special kind of maximum principle for holomorphic functions (for the proof see e.g. \cite[Lemma 4.6.6]{arebathieneu}):
\begin{lemma}\label{lem 4.6.6}
Let $U$ be an open neighborhood of $i\eta$, where $\eta\in\R$ is a given number, such that $U$ contains the closed disk $\bar{B}(i\eta ,2r)=\{ z\in\C :|z-i\eta | \le 2r\}$. Let $h:U\to \X$ be holomorphic and $c\ge 0$, $k\in \N$ such that
\begin{align*}
   \| h(z)\| \le \frac{c}{|\Re \hspace{0.3mm} z|^k}, \quad \mbox{if} \ |z-i\eta |=2r , \enskip \Re \hspace{0.3mm} z \not = 0.
\end{align*}
Then
\begin{align*}
   \| h(z)\| \le \left( \frac{4}{3}\right)^{\!k}\! \frac{c}{r^k}, \quad \mbox{for all} \enskip z\in \bar B(i\eta ,r).
\end{align*}
\end{lemma}

\section{Spectrum of a polynomially bounded sequence}

The following lemma is the key for us to set up a spectral theory for polynomially bounded sequences.
\begin{lemma} \label{lem 3.1}
Assume that $\bar S$ is the operator induced by the translation $S$ in the quotient space $l^\nu _\infty (\X) /c^\nu _0(\X)$. Then
\begin{align*}
   \sigma (\bar S) \subset \Gamma .
\end{align*}
Moreover, for each $|\lambda |\not =1$ with $|\lambda| <2$ and $f\in l^\nu_\infty (\X)$, the following estimate is valid:
\begin{align}
   \| R(\lambda , \bar S)\bar f \| \le\frac{ C}{\left||\lambda | -1 \right|^{\nu  +1}}\| f\|_\nu  ,
\end{align}
where $C$ is  a certain positive number, independent of $f$.
\end{lemma}
\begin{proof}
We will prove that if $|\lambda |\not=1$, then $\lambda \in \rho(\bar S)$. In other words, $\sigma (\bar S) \subset \Gamma $. And after that, we will give estimates of the resolvent $R(\lambda ,\bar S)\bar f$ of a given sequence $f\in l^\nu_\infty (\X)$.
To study the invertibility of the operator $(\lambda -\bar S)$, we consider the non-homogeneous linear difference equation
\begin{align}\label{3.5}
   x(n+1)-\lambda x(n) = f(n), \quad n\in \N.
\end{align}
To prove that $\lambda \in \rho (\bar S)$ for each $|\lambda|\not =1$ we will show that this equation (\ref{3.5}) has a unique solution $x\in  l_\infty^\nu (\X)$ modulo $c_0^\nu (\X)$ given $f\in  l_\infty^\nu (\X)$. 

{\it We first consider the case $|\lambda |<1$}. In this case, we will use the Variation of Constants Formula
\begin{align*}
   x(n)=\lambda^{n-1}x(1) +\sum_{k=1}^{n-1} \lambda ^{n-1-k} f(k), \quad n\in \N.
\end{align*}
Since the sequence $f$ grows polynomially, the series
$
   \sum_{k=1}^{\infty} \lambda ^{n-1-k} f(k) 
$
is absolutely convergent. Also, by $|\lambda |<1$ the sequence $\{ \lambda^{n-1}x(1) \}_{n\in\N}$ is in $c_0^\nu (\X)$. Therefore, Eq.~(\ref{3.5}) has a unique solution 
\begin{align*}
   \biggl\{x_f (n):= \sum_{k=1}^{n-1} \lambda ^{n-1-k} f(k) \biggr\}_{n\in\N}\quad \mbox{ modulo} \quad    c_0^\nu (\X).
\end{align*}
Now suppose that $g$ is any element in the class $\bar f$. We will show that $\bar x_g =\bar x_f$. Or equivalently, we have to show that whenever $h\in c_0^\nu (\X)$, the sequence 
\begin{align*}
   \{ x_h (n)\}_{n\in\N}=\biggl\{ \sum_{k=1}^{n-1} \lambda ^{n-1-k} h(k) \biggr\}_{n\in\N} \in c_0^\nu (\X).
\end{align*}
In fact, as $h\in c_0^\nu (\X)$, given $\varepsilon >0$ there exists a natural number $M$ such that for all $k\ge M$, 
\begin{align*}
   \frac{|h(k)|}{k^\nu } < \frac{1-|\lambda|}{2}\varepsilon.
\end{align*}
Therefore, for all $n\ge M+1$, 
{\allowdisplaybreaks
\begin{align*}
   \frac{  \| x_h (n) \|  }{n^\nu } 
   &\le \sum_{k=1}^{M-1} \frac{ |\lambda | ^{n-1-k} }{n^\nu }\| h(k) \|   + \sum_{k=M}^{n-1} \frac{ |\lambda | ^{n-1-k} }{n^\nu } \| h(k) \|          \\
   &\le \frac{ |\lambda|^n}{n^\nu }  \sum_{k=1}^{M-1} |\lambda |^{-1-k} \| h(k) \| +   \sum_{k=M}^{n-1} |\lambda| ^{n-1-k} \frac{ \|h(k) \|  }{k^\nu }  \\
   &\le \frac{ |\lambda|^n}{n^\nu }  \sum_{k=1}^{M-1} |\lambda |^{-1-k} \| h(k) \| + \sum_{k=M}^{n-1} |\lambda| ^{n-1-k}  \frac{1-|\lambda|}{2}\varepsilon \\
   &\le  \frac{ |\lambda|^n}{n^\nu }  \sum_{k=1}^{M-1} |\lambda |^{-1-k} \| h(k) \| +\frac{\varepsilon}{2} .
\end{align*}}%
As $M$ is a fixed natural number and $|\lambda|<1$ there exists a natural number $K\ge M+1$ such that for all $n\ge K$, 
\begin{align*}
   \frac{ |\lambda|^n}{n^\nu }  \sum_{k=1}^{M-1} |\lambda |^{-1-k} \| h(k) \|  
   \le \frac{\varepsilon}{2} .
\end{align*}
Consequently, given any $\varepsilon>0$ there exists a number $K$ such that for all $n\ge K$,
\begin{align*}
   \frac{  \| x_h (n) \|  }{n^\nu } \le \varepsilon .
\end{align*}
This means
\begin{align*}
   \lim_{n\to\infty} \frac{  \| x_h (n) \|  }{n^\nu }  =0.
\end{align*}
By this we have proved that  $\bar x_f =\bar x_g$ whenever $\bar f=\bar g$. Namely, we have showed that if $|\lambda|<1$, then $\bar x_f =(\lambda -\bar S)^{-1}f$.  In other words, $\lambda \in \rho (\bar S)$. Moreover, for any representative $g$ of the class $\bar f$
\begin{align*}
   \| R(\lambda , \bar S) \bar f \|_\nu  
    = \| \bar x_f\|_\nu  
    =\inf_{g\in \bar f} \| x_g\| _\nu 
   &\le \| x_g\| _\nu 
    = \sup_{n\in\N} \frac{\| \sum_{k=1}^{n-1} \lambda ^{n-1-k} g(k)\|}{n^\nu } \nonumber \\
   &\le \sup_{n\in\N}  \sum_{k=1}^{n-1} |\lambda| ^{n-1-k} \frac{\| g(k)\|}{k^\nu } \nonumber \\
   &\le \biggl( \sup_{n\in\N}  \sum_{k=1}^{n-1} |\lambda |^{n-1-k} \biggr) \| g\|_\nu \\
   &\le \frac{ \| g\|_\nu  }{1-|\lambda|} .
\end{align*}
Finally, as $g$ is any representative of the class $\bar f$, we have
\begin{align*}
   \| R(\lambda , \bar S) \bar f \| 
   \le \inf_{g\in \bar f} \frac{ \| g\|_\nu  }{1-|\lambda|} 
   = \frac{ \| \bar f \|_\nu  }{1-|\lambda|} .
\end{align*}

{\it Next, we consider the case $|\lambda |>1$}. We can verify that the formula
\begin{align}\label{3.12}
   x(n)= \lambda^{n-1}x(1) -\sum_{k=n}^\infty \lambda^{n-k-1}f(k), \quad n\in \N,
\end{align}
gives the general solution to Eq.~(\ref{3.5}). In fact, since $|\lambda|>1$ and $f$ grows polynomially the series $\sum_{k=n}^\infty \lambda^{n-k-1}f(k)$ is absolutely convergent for each $n\in\N$. Moreover, by (\ref{3.12}), for each $n\in \N$,
{\allowdisplaybreaks
\begin{align*}
   x(n+1) 
   &= \lambda^nx(1)-\sum_{k=n+1}^\infty \lambda^{n-k}f(k) \\
   &= \lambda^nx(1)- \sum_{k=n}^\infty \lambda^{n-k}f(k) +f(n)\\
   &= \lambda x(n)+ f(n).
\end{align*}}%
Given $f\in l_\infty^\nu (\X)$, the only solution of Eq.~(\ref{3.5}) in $l_\infty^\nu (\X)$ is
\begin{align*}
   x_f:=\biggl\{ -\sum_{k=n}^\infty \lambda^{n-k-1}f(k)\biggr\} _{n\in\N} .
\end{align*}
Indeed,
{\allowdisplaybreaks
\begin{align}
   \| x_f\|_\nu  
   &\le \sup_{n\in\N} \frac{  \sum_{k=n}^\infty |\lambda|^{n-k-1}\|f(k)\| }{n^\nu  }
    = \sup_{n\in\N}   \sum_{k=n}^\infty   \frac{    |\lambda|^{n-k-1}k^\nu }{n^\nu  }  \frac{\|f(k)\| }{k^\nu } \nonumber \\
   &\le \sup_{n\in\N}  \sum_{k=n}^\infty |\lambda|^{n-k-1}  \frac{ k^\nu }{n^\nu  } \| f\|_\nu  
    = \sup_{n\in\N}  \sum_{j=1}^\infty |\lambda|^{-j}  \left( 1+\frac{j-1}{n}\right)^\nu  \| f\|_\nu  \nonumber \\
   &\le  \sum_{j=1}^\infty |\lambda|^{-j}  j^\nu   \| f\|_\nu  .\label{3.13}
\end{align}}%
We are interested in the behavior of 
$
   \sum_{j=1}^\infty |\lambda|^{-j}  j^\nu  
$
as $|\lambda|$ gets closer and closer to $1$ (and $\infty$, respectively). To this end, we note that for each $j\in \N$,
\begin{align*}
   |\lambda|^{-j-1}  j^\nu  &\le \int_{j}^{j+1} |\lambda|^{-t}t^\nu dt.
\end{align*}
Therefore,
\begin{align*}
   \sum_{j=1}^\infty |\lambda|^{-j}  j^\nu  
   \le |\lambda| \int^\infty_0|\lambda|^{-t}t^\nu  dt 
   = |\lambda| \int^\infty_0 e^{- t\cdot \ln (|\lambda | ) }t^\nu dt 
   = \frac{  |\lambda |\nu !}{ |\ln (|\lambda|)|^{\nu  +1}}.
\end{align*}
Consequently, by using the series
\begin{align*}
   \ln (1+x )
   =\sum_{n=1}^\infty (-1)^{n+1} \frac{x^n}{n}
   =x-\frac{x^2}{2}+\frac{x^3}{3}-\cdots, 
\end{align*} 
we can show that there exists a positive number $C$ independent of $f$ such that for $1< |\lambda |<2$,
\begin{align}
   \| x_f\|_\nu  
   \le \frac{  |\lambda |\nu !}{|\ln (|\lambda|) |^{\nu  +1}}\| f\|_\nu 
   \le \frac{ C}{||\lambda | -1 |^{\nu  +1}}\| f\|_\nu  .\label{3.14}
\end{align}

Similarly as in the previous case where $|\lambda|<1$, we will prove that $\bar x_f=\bar x_g$ whenever $\bar f=\bar g$. Namely, if $\bar h=0$, then $\bar x_h=0$. In fact, for a given $\varepsilon>0$, there exists a natural number $N$ such that for all $k\ge N$, $\|h(k)\| /k^\nu  < \varepsilon$. Therefore, for all $n\ge N$,
{\allowdisplaybreaks
\begin{align*}
   \frac{\|x_h(n)\| }{n^\nu }
   &\le \sum_{k=n}^\infty \frac{ |\lambda|^{n-k-1}}{n^\nu } \| h(k)\| \\
   &\le \sum_{k=n}^\infty  \frac{  |\lambda|^{n-k-1}}{n^\nu  }\varepsilon k^\nu  \\
   &\le \varepsilon |\lambda |^{-1} \sum_{j=0}^\infty |\lambda|^{-j}\left(1+\frac{j}{n}\right)^\nu.
\end{align*}}%
Since $|\lambda |>1$ is fixed, the series $\sum_{j=0}^\infty |\lambda|^{-j}(1+j/n)^\nu $ is convergent, so this shows that
\begin{align*}
   \lim_{n\to \infty}  \frac{\|x_h(n)\| }{n^\nu } =0.
\end{align*}
That is, $\bar x_h=0$. This yields that $\lambda \in \rho (\bar S)$ and $\bar x_f =(\lambda -\bar S)^{-1}f$. Finally, with (\ref{3.14}) the proof of the lemma is complete.
\end{proof}

\begin{definition}
Let $f\in l_\infty^\nu (\X)$ be a given sequence in $\X$. Then its spectrum is defined to be the set of all complex $\xi_0\in \Gamma$ such that the complex function $R(\lambda, \bar S)\bar f$ has no analytic extension to any neighborhood of $\xi_0$. The spectrum of a sequence $f\in l_\infty^\nu (\X)$ will be denoted by $\sigma_\nu (f)$.
\end{definition}

Before we proceed we introduce some notations: $D_{|z|>1}:= \{ z\in \C : |z|>1\}$, $B(\xi_0,\delta ):=\{ z\in \C: |z-\xi_0| <\delta \}.$
\begin{lemma}
Let $f\in l_\infty^\nu (\X)$. Then, $\xi_0\in \Gamma$ is in $\sigma_\nu (f)$ if and only if the function $g : D_{|z|>1}  \ni \lambda \mapsto R(\lambda, \bar S)\bar f \in l_\infty^\nu (\X)$ cannot be extended to an analytic function in any neighborhood of $\xi_0$.
\end{lemma}
\begin{proof}
It suffices to show that if $g$ can be extended to an analytic function in a neighborhood of $\xi_0$, then $\xi_0\notin \sigma_\nu (f)$. Suppose that $g(\lambda )=h(\lambda )$ for all $\lambda \in D_{|z|>1} \cap B(\xi_0,\delta )$ where $h$ is an analytic function  in a small disk $B(\xi_0,\delta )$. Then, the function $(\lambda -\bar S)h(\lambda )$ is analytic in $B(\xi_0,\delta )$. We observe that, for $\lambda \in D_{|z|>1} \cap B(\xi_0,\delta )$
\begin{align*}
   (\lambda -\bar S)h(\lambda )=(\lambda -\bar S)g(\lambda ) =(\lambda -\bar S)R(\lambda, \bar S)\bar f =\bar f.
\end{align*}
That is, the function $(\lambda -\bar S)h(\lambda )$ is a constant in an open and connected subset $D_{|z|>1} \cap B(\xi_0,\delta )$ of the disk
$B(\xi_0,\delta )$. Hence, $(\lambda -\bar S)h(\lambda )=\bar f$ for all $\lambda$ in $B(\xi_0,\delta )$. In particular, when $|\lambda |<1$ and $\lambda \in B(\xi_0,\delta)$, $h=R(\lambda ,\bar S)\bar f$. That means, $h(\lambda )$ is an analytic extension of the function $R(\lambda ,\bar S)\bar f$ as a complex function on $\{z\in \C: |z|\not= 1\}$ to a neighborhood of $\xi_0$.
\end{proof}

\begin{proposition}\label{pro 1}
Let $f\in l_\infty^\nu (\X)$ be a given sequence in $\X$. Then the following assertions are valid:
\begin{enumerate}
   \item $\sigma_\nu (f)$ is a closed subset of $\Gamma;$
   \item The sequence $f$ is in $c_0^\nu(\X)$ if and only if $\sigma_\nu (f)=\emptyset;$
   \item If $\xi_0$ is an isolated element of $\sigma_\nu (f)$, then the point $\xi_0$ is a pole of the complex function $R(\lambda, \bar S)\bar f$ of order up to $\nu+1$.
\end{enumerate}
\end{proposition}
\begin{proof}
Part (i) is obvious from the definition of the spectrum of $x$.

Part (ii): Clearly, if $f\in c_0^\nu (\X)$, then $\sigma_\nu (f)=\emptyset$. Conversely, if  $\sigma_\nu (f)=\emptyset$, then the complex function $\hat f(\lambda ):= R(\lambda, \bar S)\bar f$ is an entire function. Moreover, it is bounded. In fact, from (\ref{3.13}) for large $|\lambda |>2$, 
{\allowdisplaybreaks
\begin{align*}
   \| \hat f(\lambda ) \|_\nu 
   \le \| \bar x_f\| 
   &\le  \sum_{j=1}^\infty |\lambda|^{-j}  j^\nu   \| f\|_\nu \\
   &= |\lambda |^{-1}  \sum_{k=0}^\infty |\lambda|^{-k}(k+1)^\nu   \| f\|_\nu \\
   &\le |\lambda |^{-1}  \sum_{k=0}^\infty \frac{(k+1)^\nu }{2^{k} }  \| f\|_\nu.
\end{align*}}%
Since the series $ \sum_{k=0}^\infty (k+1)^\nu  /2^k $ is convergent, we have
\begin{align*}
   \lim_{|\lambda | \to \infty} \| \hat f(\lambda ) \|_\nu =0.
\end{align*}
By the Liouville Theorem, this complex function $\hat f(\lambda ):= R(\lambda, \bar S)\bar f$ is the zero function, so $\bar f=0$ since $R(\lambda, \bar S)$ is injective for each large $|\lambda |$. That means $f\in c_0^\nu(\X)$.

Part (iii): Without loss of generality we may assume that $\xi_0=1$. Consider $\lambda $ in a
small neighborhood of $1$ in the complex plane. We will express
$\lambda =e^{z}$ with $|z|< \delta_0 $. Choose a small $\delta_0
>0$ such that if $|z|<\delta_0$, then
\begin{align*}
   \frac{1}{|1-|\lambda ||} \le \frac{2}{|\Re \hspace{0.3mm} z|}.
\end{align*}
It follows from Lemma \ref{lem 3.1} that for $0< |\Re \hspace{0.3mm} z |<\delta _0$, 
\begin{align*}
   \| R(\lambda ,\bar S)\bar x \| 
   \le \frac{C}{|1-|\lambda ||^{\nu+1} } \| \bar x\| 
   \le \frac{C2^{\nu+1} }{|\Re \hspace{0.3mm} z| ^{\nu +1}}  \| \bar x\| .
\end{align*}
Set $f(z )= R(e^z ,\bar S)\bar x $ with $|z|<\delta_0$. Since $1$ is
a singular point of $\| R(\lambda ,\bar S)\bar x \|$, $0$ is a
singular point of $f(z)$ in $\{ |z|<\delta_0 \}$ . For each $n\in\Z$ and $0<r<\delta_0$, we have
\begin{align*}
   \left\| \frac{1}{2\pi i} \int_{|z |=r}  \left( 1+\frac{ z  ^2}{r^2}\right)^{\nu +1} \! f(z)dz \right\| 
   \le \frac{1}{2\pi} \int_{|z |=r}  \left| 1+\frac{ z  ^2}{r^2}\right|^{\nu +1} \| f(z)\| \, |dz | .
\end{align*}
If $z =re^{i\varphi}$, where $\varphi$ is real, one has
\begin{align*}
   \left| 1+\frac{ z ^2}{r^2}\right| ^{\nu+1}\!\!
   &= |1+e^{2i\varphi} |^{\nu+1}
    = | e^{-i\varphi}+e^{i\varphi}|^{\nu+1} \\
   &= (2| \cos \varphi|)^{\nu+1} 
    = 2^{\nu+1}r^{ -\nu-1}|\Re \hspace{0.3mm} z|^{\nu+1}.
\end{align*}
Therefore,
\begin{align}
   &\left\| \frac{1}{2\pi i} \int_{|z |=r}  \left( 1+\frac{z^2}{r^2}\right)^{\nu +1} \! \frac{f(z)}{z^{n+1}}dz \right\| \nonumber \\
   &\le \frac{1}{2\pi} \int_{|z|=r}   2^{\nu+1}r^{ -n-\nu-2}|\Re \hspace{0.3mm} z|^{\nu+1}    \frac{C2^{\nu+1} }{|\Re \hspace{0.3mm} z| ^{\nu +1}}  \| \bar x\| \, |dz | \nonumber \\
   &= \frac{ C4^{\nu+1} r^{-n-\nu-2}}{2\pi} \int_{|z|=r}  |dz | \, \| \bar x\| \nonumber \\
   &= C4^{\nu+1}r^{-n-\nu-1} \| \bar x\| . \label{15}
\end{align}
Consider the Laurent series of $f(z )$ at $z =0$,
\begin{align*}
   f(z) = \sum_{n=-\infty}^{\infty} a_n z^n,
\end{align*}
where
\begin{align*}
   a_n =\frac{1}{2\pi i} \int_{| z |=r} \frac{f(z) }{z ^{n+1}}dz, \quad n\in \Z .
\end{align*}
It follows that for each $n\in \Z$,
\begin{align*}
   \frac{1}{2\pi i} \int_{|z |=r}  \left( 1+\frac{z^2}{r^2}\right)^{\nu +1} \frac{f(z)}{z^{n+1}}dz 
   &= \frac{1}{2\pi i} \int_{|z |=r}   \sum_{k=0}^{\nu+1} \frac{(\nu+1)!}{k!(\nu+1-k)!} r^{-2k}     \frac{f(z)}{z^{n+1-2k }}  dz  \\
   &= \sum_{k=0}^{\nu+1} \frac{(\nu+1)!}{k!(\nu+1-k)!} r^{-2k} \frac{1}{2\pi i} \int_{|z |=r}   \frac{f(z)}{z^{n+1-2k }} dz  \\
   &= \sum_{k=0}^{\nu+1} \frac{(\nu+1)!}{k!(\nu+1-k)!} r^{-2k} a_{n-2k}.
\end{align*}
This, together with (\ref{15}), shows
\begin{align*}
   \left|\,\sum_{k=0}^{\nu+1} \frac{(\nu+1)!}{k!(\nu+1-k)!} r^{-2k} a_{n-2k} \right| 
   \le C4^{\nu+1}r^{-n-\nu-1}\| \bar x\| .
\end{align*}
Multiplying both sides by $r^{2\nu}$ gives
\begin{align*}
   \left|\,\sum_{k=0}^{\nu+1} \frac{(\nu+1)!}{k!(\nu+1-k)!} r^{2\nu-2k} a_{n-2k} \right| 
   \le C4^{\nu+1}r^{\nu-n-1}\| \bar x\|.
\end{align*}
Observe that in the left side is a polynomial in terms of $r$ whose zero power term is $a_{n-2\nu}$. Therefore, when $\nu -n-1\ge 1$ if we let $r$ to get closer and closer to zero, then $a_{n-2\nu}$ must be zero. That is for all $\nu \ge n+2$, the coefficients $a_{n-2\nu}=0$. This yields that for all $j\le - \nu -2$, $a_j=0$. In other words, $z=0$, or $\lambda =1$ is a pole of the complex function $\hat f(\lambda ):= R(\lambda, \bar S)\bar f$ with order up to $\nu +1$.
\end{proof}

Before proceeding we introduce a notation:  Let $z\in \C\backslash\{0\}$ such that
$z= re^{i\varphi}$ with reals $r,\varphi$, and $F(z)$ be any complex
function. Then we define
\begin{align}
   \lim_{\lambda \downarrow z} F(\lambda ):= \lim_{s\downarrow r}F(s e^{i\varphi}).
\end{align}

\begin{corollary}\label{cor 4}
Let $f\in l_\infty ^\nu (\X)$, and $\xi_0\in \Gamma$ be an isolated point in $\sigma_\nu (f)$. Moreover, assume that  
\begin{align}\label{con 3.26}
   \lim_{\lambda \downarrow \xi_0} (\lambda -\xi_0)R(\lambda , \bar S)\bar f=0.
\end{align}
Then the singular point $\xi_0$ of $R(\lambda , \bar S)\bar f$ is removable and the complex function $R(\lambda , \bar S)\bar f$ is zero in the connected open subset of its domain that contains $\xi_0$.
\end{corollary}
\begin{proof}
By Proposition \ref{pro 1}, $\xi_0$ is a pole of order up to $\nu+1$. Consider the Laurent series of $R(\lambda , \bar S)\bar f$ in a neighborhood of $\xi_0$ we have
\begin{align*}
   R(\lambda , \bar S)\bar f = \sum_{j=-\nu-1}^\infty \frac{a_j}{(\lambda -\xi_0)^{j+1}}.
\end{align*} 
If (\ref{con 3.26}) is satisfied, then for any $k\ge 1$ the following is also valid:
\begin{align*}
   \lim_{\lambda \downarrow z} (\lambda -\xi_0)^kR(\lambda , \bar S)\bar f=0.
\end{align*}
If we let $k$ take on the values $1,2,\ldots ,\nu+1$, then we see that $a_j=0$ for all $j=-\nu-1,-\nu, \ldots. $ That is, the function $R(\lambda , \bar S)\bar f$ is zero (so analytic) in a neighborhood of $\xi_0$. From the properties of analytic functions this function must be zero in the connected open subset of its domain as well.
\end{proof}

\section{Applications to study the asymptotic behaviors of difference equations in Banach spaces}
\subsection{Asymptotic behavior of polynomially bounded solutions of difference equations}
In this subsection we will apply the results obtained in the previous section to study the polynomially bounded solutions of difference equations of the form
\begin{align}\label{eq1}
   x(n+1)=Tx(n) + F(n), \quad n\in \N,
\end{align}
where $T$ is a bounded linear operator in $\X$ and $F\in c_0^\nu(\X)$.

\begin{definition}
A bounded operator $T$ from a Banach space $\X$ to itself is said to be $\nu$-polynomially power bounded, if
\begin{align*}
   \sup_{n\in \N} \frac{\| T^n\|}{n^\nu } < \infty ,
\end{align*} 
where $\nu$ is a nonnegative integer.
\end{definition}

\begin{lemma}\label{lem 4.2}
Let $x\in l^\nu _\infty (\X)$ be a solution of (\ref{eq1}). Then
\begin{align}\label{spec est}
   \sigma_\nu (x) \subset \sigma (T) \cap \Gamma .
\end{align}
Moreover, for $\lambda \in \rho (\bar S)\cap \rho (T)$,
\begin{align}\label{resolvent}
   R(\lambda ,\bar S)\bar x = R(\lambda ,\bar T)\bar x.
\end{align}
\end{lemma}
\begin{proof}
Consider the operator of multiplication by $T$ in the spaces $l^\nu_\infty (\X)$. It is easy to see that the operator is bounded and preserves $c_0^\nu(\X)$, so it induces an operator $\bar T$ in the quotient space $l^\nu _\infty (\X)/c_0^\nu(\X)$. Moreover,
$
   \sigma (\bar T) \subset \sigma (T).
$
Since $x$ is a solution of (\ref{eq1}), for each $|\lambda|\not= 1$ we have 
\begin{align*}
   R(\lambda , \bar S)\bar S \bar x 
   = R(\lambda , \bar S)\bar T \bar x +R(\lambda , \bar S)\bar F
   = \bar T R(\lambda , \bar S)\bar x,
\end{align*} 
This, together with the identity $\lambda R(\lambda , \bar S)\bar x-\bar x= R(\lambda , \bar S)\bar S \bar x$, shows
\begin{align*}
   \lambda R(\lambda , \bar S)\bar x-\bar x
   = \bar T R(\lambda , \bar S)\bar x.
\end{align*} 
Therefore,
\begin{align*} 
   \bar x
   = \lambda R(\lambda , \bar S)\bar x-\bar T R(\lambda , \bar S)\bar x
   = (\lambda -\bar T) R(\lambda , \bar S)\bar x .
\end{align*} 
If $\xi_0\in \Gamma$ and $\xi_0\not\in \sigma (T)$, then there exists a neighborhood of $\xi_0$ (in $\C$) such that for any $\lambda \in U$ and $ |\lambda |\not=1$,
\begin{align}\label{resolvent2}
   R(\lambda ,\bar T)\bar x = R(\lambda ,\bar S)\bar x .
\end{align}
As the left hand side function is an analytic extension in a neighborhood $U$ of $\xi_0$, by (\ref{resolvent2}), the complex function $R(\lambda , \bar S)\bar x$ has an analytic extension to the neighborhood $U$ of $\xi_0$, that is $\xi_0\not\in \sigma_\nu (x)$. In other words, $\sigma _\nu (x) \subset \sigma (T) \cap \Gamma$. Moreover, (\ref{resolvent}) is proved.
\end{proof}

We will prove the following that extends the famous Katznelson-Tzafriri to the case of $\nu$-polynomially bounded operator.

\begin{theorem}\label{the kat ext}
Let $T\in L(\X)$ be $\nu$-polynomially bounded such that $\sigma (T) \cap \Gamma \subset \{ 1\}$, where $\nu$ is a nonnegative integer. Then 
\begin{align}\label{kat}
   \lim_{n\to \infty} \frac{1}{n^\nu}  (T-I)^{\nu +1} T^n  =0.
\end{align}
\end{theorem}
\begin{proof} We consider the sequence $x:= \{ x(n):= T^n\}_{n=1}^\infty$ in $L(\X)$. Obviously, $\{x(n)\}_{n=1}^\infty \in l^\nu_\infty (L(\X))$. Let us denote by $\bar T$ the operator of multiplication by $T$ in $\Y:=L(\X)$. Then, we have an equation in $\Y$:
\begin{align*}
x(n+1)=\bar T x(n), \quad n\in \N.
\end{align*}
Note that $\sigma (\bar T)\subset \sigma (T)$. Therefore, by Lemma \ref{lem 4.2}, we have  $\sigma_\nu (x)\subset \{1\}$. For each $\lambda \in \rho (S)$, we have the identity
\begin{align*}
   R(\lambda ,\bar S)\bar S\bar x
   = \lambda R(\lambda ,\bar S)\bar x - \bar x.
\end{align*}
By a simple induction we can show that for each $j\in\N$,
\begin{align*}
   R(\lambda , \bar S) \bar S^j \bar x
   = \lambda^j R(\lambda ,\bar S)\bar x- P(\lambda , \bar x, \bar S\bar x),
\end{align*}
where $P(\lambda , \bar x, \bar S\bar x)$ is a polynomial of $\lambda,\bar x,\bar S\bar x$. Hence,
\begin{align*}
   R(\lambda ,\bar S) (\bar S-I)^{\nu+1} \bar x
   = (\lambda -1)^{\nu+1} R(\lambda,\bar S)\bar x +Q(\lambda ,\bar x,\bar S\bar x),
\end{align*}
where $Q(\lambda ,\bar x,\bar S\bar x)$ is a polynomial of $\lambda,\bar x,\bar S\bar x$.
Note that $\sigma_\nu ((S-I)^\nu x) \subset \sigma_\nu (x) \subset \{ 1\}$. By Proposition \ref{pro 1}, $1$ is a pole of $\nu+1$ order of the complex function $g(\lambda ):=R(\lambda,\bar S)\bar x$, so the complex function
\begin{align*}
   \lambda \mapsto R(\lambda ,\bar S) (\bar S-I)^{\nu+1} \bar x
\end{align*} 
is extendable analytically to a neighborhood of $1$. Therefore, for the sequence $y:= (S-I)^{\nu+1} x$ we have $\sigma_\nu (y) =\emptyset$. By Proposition \ref{pro 1}, $y= (S-I)^{\nu+1} x \in c^\nu _0(L(\X))$, that is, (\ref{kat}) is valid.
\end{proof}

\begin{remark}
A famous Katznelson-Tzafriri Theorem (see \cite{kattza}) is stated as follows: Let $T\in L(\X)$ satisfy $\sup_{n\in \N} \|T^n\| <\infty$ and $(\sigma (T) \cap \Gamma) \subset \{1\}$. Then
\begin{align*}
   \lim_{n\to\infty} (T^{n+1}-T^n)=0.
\end{align*}
There are many extensions of this theorem (see e.g.~\cite{vu1} and its references). An elementary proof of this theorem is given in \cite{min}. 
In Theorem \ref{the kat ext}, when $\nu =0$, we obtain the above mentioned Katznelson-Tzafriri Theorem.
\end{remark}

Below is an individual version of Katznelson-Tzafriri Theorem for possibly non-$\nu$ polynomially bounded operator $T$.

\begin{theorem}\label{the kat ext2}
Let $T\in L(\X)$ satisfy $\sigma (T) \cap \Gamma \subset \{ 1\}$. Then, for each $x_0\in \X$,
\begin{align}\label{kat2}
   \lim_{n\to \infty} \frac{1}{n^\nu} (T-I)^{\nu +1} T^nx_0 =0,
\end{align}
provided that
\begin{align*}
   \sup_{n\in\N} \frac{\| T^nx_0\| } {n^\nu} < \infty .
\end{align*}
\end{theorem}
\begin{proof}
The proof is similar to that of Theorem \ref{the kat ext}. In analogy to the sequence $x$ in the proof of Theorem \ref{the kat ext} we can use the sequence $\{ T^nx_0\}_{n=1}^\infty$.
\end{proof}

Let us consider homogeneous linear difference equations of the form
\begin{align}\label{eq2}
   x(n+1)=Tx(n), \quad n\in \N,
\end{align}
where $T\in L(\X)$. Each solution $\{x(n)\}_{n=1}^\infty$ of this equation (\ref{eq2}) is of the form $x(n)=T^{n-1}x_0$, $n\in \N$, for some $x_0\in \X$. 

\begin{theorem}\label{the 1}
Let $T\in L(\X)$, and let $x$ be a $\nu$ polynomially bounded solution of Eq.~(\ref{eq1}). Assume further that the following conditions are satisfied:
\begin{enumerate}
   \item $\sigma (T) \cap \Gamma$ is countable 
   \item For each $\xi_0=e^{i\phi_0}\in \sigma (T)\cap \Gamma$ 
\begin{gather} 
   \{ z= re^{i\phi_0}, r>1\} \subset \rho(T); \label{4.5} \\
   \lim_{\lambda \downarrow z} (\lambda -\xi_0)R(\lambda , \bar T) \bar x=0. \label{4.6}
\end{gather}
\end{enumerate}
Then
\begin{align}\label{0}
   \lim_{n\to \infty} \frac{x(n)}{n^\nu} =0 .
\end{align}
\end{theorem}
\begin{proof}
Since $\sigma_\nu (f) \subset \sigma (T) \cap \Gamma $, if $ \sigma (T) \cap \Gamma $ is empty, then the claim of the theorem is clear. Next, if it is not, then from the countability of $\sigma_\nu (f) $ as a closed subset of $\Gamma$ there must be an isolated point, say $\xi_0$ of $\sigma_\nu (f) $. However, by condition (\ref{4.6}) and Corollary \ref{cor 4} the set of non-removable singular points of the complex function $R(\lambda , \bar S) \bar x=R(\lambda , \bar T) \bar x$ cannot have an isolated point. That means, $\sigma_\nu (x) $ must be empty set, so by Proposition \ref{pro 1}, the sequence $x=\{x(n)\}_{n=1}^\infty$ must be in $c_0^\nu(\X)$, that is (\ref{0}).
\end{proof}

The following result gives a sufficient condition for the stability of polynomially bounded solutions that is well known as Arendt-Batty-Ljubich-Vu Theorem (see \cite{arebathieneu}):

\begin{corollary}
Let $T\in L(\X)$ be $\nu$-polynomially power bounded. Assume further that 
\begin{enumerate}
   \item $\sigma (T)\cap \Gamma$ is countable,
   \item  For each $\xi_0$ of $ \sigma (T)\cap\Gamma,$ and each $x\in \X,$
\begin{align}\label{4.8}
   \lim_{\lambda \downarrow z} (\lambda -\xi_0)R(\lambda ,T) x =0.
\end{align}
\end{enumerate}
Then, for each $x\in \X,$
\begin{align} 
   \lim_{n\to \infty} \frac{T^nx}{n^\nu} =0 .
\end{align}
\end{corollary}
\begin{proof}
It is clear that $x(n)=T^nx$ is a solution of Eq.~(\ref{eq2}).  By the Spectral Radius Theorem the spectral radius $r_\sigma (T)$ of $T$ must satisfy $r_\sigma (T) \le 1$ because of the polynomial boundedness of $T$, so (\ref{4.5}) is satisfied. By Theorem \ref{the 1}, we only need to check condition (\ref{4.6}). We have
\begin{align*}
   0\le \lim_{\lambda \downarrow z} \| (\lambda -\xi_0) R(\lambda , \bar T) \bar x \|_\nu 
   &\le \lim_{\lambda \downarrow z} \sup_{n\in \N} \frac{\| (\lambda -\xi_0) R(\lambda ,  T) T^nx \|}{n^\nu} \\
   &\le  \lim_{\lambda \downarrow z} \sup_{n\in\N} \frac{\| T^n\|}{n^\nu}   \| (\lambda -\xi_0)    R(\lambda ,  T) x \| \\
   &=  \sup_{n\in\N} \frac{\| T^n\|}{n^\nu}   \lim_{\lambda \downarrow z}   \| (\lambda -\xi_0) R(\lambda ,  T) x \|.
\end{align*}
Since $T$ is $\nu$-polynomially power bounded $\sup_{n\in\N} (\| T^n\|/n^\nu )$ is finite, so (\ref{4.8}) yields that condition (\ref{4.6}) is satisfied.
\end{proof}

\subsection{Asymptotic behavior of solutions of fractional difference equations}

Consider fractional difference equations of the form
\begin{align}\label{frac eq}
   \Delta^\alpha x(n) =Tx(n)+y(n), \quad n\in \N,
\end{align}
where $0<\alpha\leq 1$, $T\in L(\X)$ and $y\in c_0^\nu(\X)$.

\begin{definition} \label{def 4.9}
(\cite[Definition 3.1]{liz}) Let $T$ be a bounded operator defined on a Banach space $\X$ and $\alpha>0$. We call $T$ the {\it generator of an $\alpha$-resolvent sequence} if there exists a sequence of bounded and linear operator $\{ S_\alpha (n)\}_{n\in \N} \subset L(\X)$ that satisfies the following properties
\begin{enumerate}
   \item $S_\alpha (0)=I;$
   \item $S_\alpha (n+1)= k^\alpha (n+1)I+T\sum_{j=0}^n k^\alpha (n-j)S_\alpha (j), \ \mbox{for all} \ n\in \N .$
\end{enumerate}
\end{definition}

As shown in \cite[Theorem 3.4]{liz} and a note before it, $S_\alpha$ is determined by one of the following formulas:

\begin{theorem} \label{thm 4.9}
Let $\alpha > 0$ and $T$ be a bounded operator defined on a Banach space $\X$. The following properties are equivalent:
\begin{enumerate}
   \item $T$ is the generator of an $\alpha$-resolvent sequence $\{S_\alpha (n)\}_{n\in \N};$
   \item $$S_\alpha (n) =\sum_{j=0}^n \frac{\Gamma (n-j+(j+1)\alpha)}{\Gamma (n-j+1)\Gamma (j\alpha + \alpha)} T^j ;$$
   \item $$S_\alpha (n) =\frac{1}{2\pi i} \int_{C} z^n ((z-1)^\alpha z^{1-\alpha}-T)^{-1}dz ,$$
where $C$ is a circle, centered at the origin of the complex plane, that encloses all spectral values of $(z-1)^\alpha z^{1-\alpha}-T.$
\end{enumerate}
\end{theorem}

\begin{theorem}{\rm(\cite[Theorem 3.7]{liz})} \label{frac sol rep}
Let $0<\alpha <1$ and $\{ y(n)\}_{n\in \N}$ is given. The unique solution of Eq.~(\ref{frac eq}) with initial condition $u(0)=x$ can be represented by
\begin{align*}
   u(n)=S_\alpha (n)u(0) +(S_\alpha * y)(n-1), \quad \mbox{for all} \enskip n\in \N.
\end{align*}
\end{theorem}


Recall that the $Z$-transform of a sequence $x:=\{ x(n)\}_{n=0} ^\infty$ is defined as
\begin{align}
   \tilde{x}(z):=\sum_{j=0}^\infty x(j)z^{-j}.
\end{align}
Let us denote $D_{|z|>1}:=\{ z\in \C: \, |z|> 1\}$, and $D_{|z|<1}:= \{ z\in \C : \, |z| <1\}.$ For each $\{x(n)\}_{n\in \N} \in l^\nu_\infty (\X)$ we will set $x(0)=0$, so some properties of the Z-transform of sequences can be stated in the following:

\begin{proposition}
Let $\{x(n)\}_{n\in \N}$ and $\{ y(n)\}_{n\in \N}$ be in $l_\infty ^\nu (\X)$. Then
\begin{enumerate}
   \item $\tilde{x}(z)$ is a complex function in $z\in D_{|z|>1};$
   \item $\widetilde{Sx}(z)= z\tilde{x}(z)-zx(0);$
   \item $\widetilde{x*y}(z)=\tilde{x}(z)\cdot \tilde{y}(z)$.
\end{enumerate}
\end{proposition}
\begin{proof}
For the proof see e.g.~\cite[Chapter 6]{ela}.
\end{proof}

There is an interesting relation between the Z-transform and Laplace transform of a function via the Poisson distribution as it is discussed in \cite{liz2}. As shown in \cite[Example 3.3]{liz2}, for the sequence $\{k^\alpha (n)\}_{n=1}^\infty$ (see (\ref{k alpha})) its Z-transform is given by
\begin{align}\label{z-trans of k}
   \tilde k^\alpha (z)= \frac{z^\alpha}{(z-1)^\alpha} .
\end{align}
We note that on the unit circle $\Gamma$, the complex function  $\tilde k^\alpha (z)$ is analytic everywhere except for the point $z=1$.
To study fractional difference equations (\ref{frac eq}) we will need the following analog of \cite[Lemma 3.3]{min2}:

\begin{lemma}\label{lem 4.13}
Let $\{ x(n)\}_{n\in \N} \in l_\infty ^\nu (\X)$. If the Z-transform $\tilde{x}(z)$ of the sequence $x$ has a holomorphic extension to a neighborhood of $z_0\in\Gamma$, then $z_0\not\in \sigma _\nu (x).$
\end{lemma}
\begin{proof}
Assume that $\tilde{x}(z)$ (with $|z|> 1$) can be extended to a holomorphic function $g_0(z)$ in $B(z_0,\delta )$ with a sufficiently small positive $\delta$. We will show that $R(z,S)x$ (with $|z|>1$) has a holomorphic extension in a neighborhood of $z_0$. By setting $x(0)=0$ we define a sequence $\{ g_k(z)\}_{k=1}^\infty $ as follows:  
\begin{align}
   g_k(z):= z^{k-1}\tilde x (z) -\sum_{j=0}^{k-1} z^{k-1-j}x(j), \quad k\in \N .
\end{align}
We are going to prove that this defines a bounded function $g(z)$ with $z$ in a small disk $B(z_0,\delta ):=\{ z\in \C: \, |z-z_0| < \delta\} $, and then, applying a necessary and sufficient condition for a locally bounded function to be holomorphic to prove that $R(z,S)x$ is holomorphic. To prove the boundedness of $g(z)$ in a small disk  $B(z_0,\delta )$ we will use a special maximum principle as in \cite{min2}. We have
\begin{align*}
   R(z,S)x 
   &:= (z-S)^{-1}x
     =  z^{-1}(I-z^{-1}S)^{-1}x \\
  &\ = z^{-1}\sum_{n=0}^{\infty}z^{-n}S(n)x
     = \sum_{n=0}^{\infty}z^{-n-1}S(n)x .
\end{align*}
Therefore, for $z\in B(z_0,\delta)\cap D_{|z|>1}$ and for each $k\in \N$,
\begin{align*}
   [R(z,S)x] (k) 
   = \sum_{n=0}^{\infty}z^{-n-1}x(n+k) 
   =  z^{-1}  \biggl( z^k\tilde {x}(z)-\sum_{j=0}^{k-1} z^{k-j}x(j) \biggr) 
   = g_k(z) .
\end{align*}
By (\ref{3.13}) and (\ref{3.14}), for $z\in B(z_0,\delta)\cap D_{|z|>1}$, there is a certain number $C$ such that
\begin{align}
   \sup_{k\in\N} \frac{\| g_k(z)\|}{k^\nu } 
    = \| g(z) \|_\nu 
   &= \biggl\| \biggl\{  \sum_{n=0}^{\infty}z^{-n-1}x(n+k)\biggr\}_{k=1}^\infty \biggr\|_\nu \nonumber \\
   &\le \frac{ C}{(|z | -1 )^{\nu  +1}}\| z\|_\nu  .\label{4.25}
\end{align}

On the other hand, for $z\in B(z_0,\delta)\cap D_{|z|<1}$ we have for all $k\in \N$,
{\allowdisplaybreaks
\begin{align*}
   \| g_k(z) \| 
   &\le |z|^{k-1} \| g_0(z)\| + \sum_{j=0}^{k-1} |z|^{k-1-j} \| x(j)\| \\
   &\le |z|^{k-1} \| g_0(z)\| +\sum_{j=0}^{k-1} |z|^{k-1-j}j^\nu \| x\|_\nu \\
   &\le \biggl( \sup_{z\in B(z_0,\delta )} \| g_0(z)\| +\| x\|_\nu   \biggr) \sum_{j=0}^{k-1}|z|^{k-1-j}j^\nu \\
   &= M \sum_{j=0}^{k-1}|z|^{k-1-j}j^\nu ,
\end{align*}}%
where $ M:= \sup_{z\in B(z_0,\delta )} \| g_0(z)\| +\| x\|_\nu.$
Hence, for all $k\in \N$,
\begin{align}
   \frac{\| g_k(z) \|}{k^\nu}  
   \le M \sum_{j=0}^{k-1}|z|^{k-1-j}\left(\frac{j}{k} \right)^\nu 
   \le M \sum_{j=0}^{k-1}|z|^{k-1-j} 
   \le \frac{M}{1-|z|}.\label{4.26}
\end{align}
By (\ref{4.25}) and (\ref{4.26}) we have proved that there is a positive number $K$ such that for $z\in B(z_0,\delta)$ and and for each $k\in \N$, this estimate is valid: 
\begin{align}
   \frac{\|g_k(z)\|}{k^\nu} \le  \frac{ K}{||z | -1 |^{\nu  +1}} .
\end{align}
Applying the maximum principle Lemma \ref{lem 4.6.6} as in \cite{min2} to the function $g_k(z)/k^\nu$ gives the boundedness of $g_k(z)/k^\nu$ in $B(z_0,\delta/2)$. In fact,
it is clear that for each $k\in \N$ the function $g_k(z)/k^\nu$ is holomorphic in $z\in B(z_0,\delta)$. Therefore, $g_k(z)/k^\nu$ is bounded by a number independent of $k$, so $g(z)$ is bounded in $B(z_0,\delta/2)$. We are now ready to apply Theorem \ref{the anal}, a criterion for a locally bounded function to be holomorphic. In fact, since the family $W:= \{ x^*\circ p_k, \ x^*\in \X^*, \ p_k: \{ x_n\}\mapsto x_k, \,k\in \N\}$ is separating and $x^*\circ p_k (g(\cdot )) = x^* (g_k(\cdot ))$ is holomorphic, the complex function $g(z)$ is holomorphic for $z\in B(z_0,\delta/2)$.

At this point we have shown that $g(z)$ is holomorphic for $z\in B(z_0,\delta/2)$, and as $g(z)= R(z,S)x$ for $|z|>1$. This yields that $R(z,S)x$ has a holomorphic extension $g(z)$ to a neighborhood of $z_0$. This completes the proof of the lemma.
\end{proof}

\begin{definition}
We denote by $\sigma_{Z,\nu} ( x)$ the set of all points $\xi_0$ on $\Gamma$ such that the Z-transform of a sequence $x:=\{x(n)\}_{n\in \N} \in l^\infty _\nu (\X)$ cannot be extended holomorphically to any neighborhood of $\xi_0$, and call this set the Z-spectrum of the sequence $x$.
\end{definition}

In the simplest case where $\nu=0$, $\sigma_\nu (x)$ may be different from $\sigma_{Z,\nu} ( x)$. In fact, the following numerical sequence
$x:=\{x(n)\}_{n\in \N} \in l^\infty_0(\R)$, where
\begin{align*}
   x(n):= 
   \begin{cases} 
      0, &n=0, \\
      1/n, &n\in \N 
   \end{cases}
\end{align*}
is in $c_0(\R)$. Obviously, $\bar x =0$, so $\sigma (x)=\emptyset$. However, $ 1 \in \sigma_{Z,\nu}(x)$ because $\tilde{x}(z) =\sum_{j=1}^\infty z^{-j} /j$ cannot be extended holomorphically to a neighborhood of $1$. In general, we only have the following inclusion.

\begin{corollary}\label{cor 4.15}
For each $x:=\{x(n)\}_{n\in \N} \in l^\infty _\nu (\X),$
\begin{align*}
   \sigma _\nu (x) \subset \sigma _{Z,\nu} (x) .
\end{align*}
\end{corollary}
\begin{proof}
The corollary is an immediate consequence of Lemma \ref{lem 4.13} and the definitions of the spectra mentioned in the statement.
\end{proof}

Before we proceed, we introduce a notation 
\begin{align*}
   \Sigma_0 := 
   &\{ z_0\in \Gamma \subset \C: \ \mbox{$(z-\tilde k^\alpha (z)T)^{-1}$ exists and $(z-\tilde k^\alpha (z)T)^{-1}$ and $\tilde k^\alpha (z)$} \\   
   &\ \mbox{are holomorphic in a neighborhood of $z_0$} \} 
\end{align*}
and $\Sigma =\Gamma \backslash \Sigma_0$. By (\ref{z-trans of k}), $\Sigma_0$ can be defined as 
\begin{align*}
   \Sigma_0 := 
   &\{ z_0\in \Gamma \subset \C: \ \mbox{$(z^{1-\alpha}(z-1)^\alpha -T)^{-1}$ \ exists and} \\   
   &\ \mbox{is holomorphic in a neighborhood of $z_0$} \} 
\end{align*}

\begin{lemma}
Let $\alpha >0$ and $S_\alpha := \{ S_\alpha (n)\}_{n\in \N} \subset L(\X)$ be the resolvent of Eq.~(\ref{frac eq}) that satisfies 
\begin{align*}
   \sup_{n\in \N} \frac{\| S_\alpha (n)\|}{n^\nu} < \infty .
\end{align*}
Then
\begin{align*}
   \sigma _\nu (S_\alpha ) \subset \Sigma .
\end{align*}
\end{lemma}
\begin{proof}
It suffices to show that if $z_0\in \Sigma_0$, then $z_0\not\in \sigma_\nu (S_\alpha)$. Taking the Z-transform of $S_\alpha$ from the equation in Definition \ref{def 4.9} gives
\begin{align*}
   z\tilde S_\alpha (z)-zS_\alpha (0) 
   = \widetilde{S S_\alpha}(z)
   = z\tilde k^\alpha (z) I -zk^\alpha (0)I + \tilde k ^\alpha (z) \cdot T\tilde S_\alpha (z) .
\end{align*}
Therefore, for $z\in D_{|z|>1}$,
\begin{align*}
   (z-\tilde k^\alpha (z)T) \tilde S_\alpha (z)
   = zS_\alpha (0) +z\tilde k^\alpha (z) I -zk^\alpha (0)I .
\end{align*}
Let $z_0\in \Sigma_0$. Then $(z-\tilde k^\alpha (z)T)^{-1}$ exists. Hence,
\begin{align*}
   \tilde S_\alpha (z)
   = (z-\tilde k^\alpha (z)T)^{-1}(zS_\alpha (0) +z\tilde k^\alpha (z) I -zk^\alpha (0)I ).
\end{align*}
And, it is clear that $ \tilde S_\alpha (z)$ has a holomorphic extension to a neighborhood of $z_0$ because both $\tilde k^\alpha (z)$ and $(z-\tilde k^\alpha (z)T)^{-1}$ are holomorphic in a neighborhood of $z_0$, so $z_0\not\in \sigma _{Z,\nu }(S_\alpha)$. By Corollary \ref{cor 4.15} this yields $z_0\not\in \sigma_\nu (S_\alpha)$. This completes the proof of the lemma.
\end{proof}

\begin{theorem}\label{the kat for frac}
Let $0<\alpha \le 1$ and $\Sigma \subset \{ 1\}$. Assume further that the $\alpha$-resolvent $S_\alpha$ of Eq.~(\ref{frac eq}) satisfies
\begin{align*}
   \sup_{n\in \N} \frac{\| S_\alpha (n)\|}{n^\nu} < \infty .
\end{align*}
Then
\begin{align}\label{4.32}
   \lim_{n\to \infty} \frac{1}{n^\nu} \sum_{k=0}^{\nu+1} \frac{(\nu+1)!}{k!(\nu+1-k)!} (-1)^{\nu+1+k} S_\alpha (n+k) =0.
\end{align}
\end{theorem}
\begin{proof}
As in the proof of Theorem \ref{the kat ext}, we can show that
\begin{align*}
   \lambda \mapsto R(\lambda ,\bar S) (\bar S-I)^{\nu+1} \bar S_\alpha
\end{align*}
has a holomorphic extension to a neighborhood of $1$ in the complex plane. Moreover, since $\sigma_\nu (S_\alpha )\subset \Sigma$ this function has a holomorphic extension to a neighborhood of all points of $\Gamma$. Namely, $\sigma_\nu ((\bar S-I)^{\nu+1} \bar S_\alpha )=\emptyset$. Therefore, $(S-I)^{\nu+1} S_\alpha \in c_0^\nu (\X)$. In other words,
\begin{align}
   \lim_{n\to \infty} \frac{1}{n^\nu} \left(  (S-I)^{\nu+1} S_\alpha  \right)(n) =0. \label{main}
\end{align}
Then it follows that 
\begin{align*}
   (S-I)^{\nu+1}
   = (-1)^{\nu+1}(I-S)^{\nu+1} 
   = (-1)^{\nu+1}\sum_{k=0}^{\nu+1} \frac{(\nu+1)!}{k!(\nu+1-k)!}(-1)^{k}S^k.
\end{align*}
This, together with (\ref{main}), yields (\ref{4.32}). The theorem is proved.
\end{proof}

\begin{remark}
When $\alpha =1$ Eq.~(\ref{frac eq}) becomes
\begin{align*}
   x(n+1)=(I+T)x(n)+y(n), \quad n\in \N.
\end{align*}
As shown in \cite{liz} $S_\alpha (n) =(I+T)^n$, $n\in\N$. With this formula, (\ref{4.32}) becomes
\begin{align*}
   \lim_{n\to \infty}\frac{1}{n^\nu}T^{\nu+1}(I+T)^n=0.
\end{align*}
Hence, Theorem \ref{the kat for frac} coincides with Theorem \ref{the kat ext} when $\alpha =1$. In other words, Theorem \ref{the kat for frac} is an extension of the Katznelson-Tzafriri Theorem for fractional difference equations (\ref{frac eq}).
\end{remark}

\section{Examples}

In this section we will give some examples to illustrate our obtained results. As is well known, the asymptotic behavior of solutions to the equation 
\begin{equation}\label{1}
   x(n+1)=Ax(n), \quad n\in \N,
\end{equation}
where $A$ is a matrix, $x\in \C^d$, depends on $\sigma (A)$ (in this case it is the set of the eigenvalues of $A$). The solutions are uniformly exponentially stable (that is, there exist positive constants $N,$ $\alpha$ such that $\|A^nx\| \le N e^{-\alpha n}\|x\|$ for all $n\in\N$, $x\in \C^d$) if $\sigma (A) \subset \{ z\in \C : |z|<1\}$. This result for the general case of Banach space $\X$ can be proved easily by the Spectral Radius Theorem as well. The asymptotic behavior of solutions of Eq.~(\ref{1}) is complicated and interesting when there are some eigenvalues lying on the unit circle. 
For example, in \cite[Theorem 1.100]{goopet} the zero solution of Eq.(\ref{1}) is stable (that is there is a constant $N$ such that $\| A^nx\| \le N\| x\|$ for all $x\in \C^d$, $n\in\N$)  when the matrix $A$ has all eigenvalues $\lambda_k$, $k=1,2,\dots,n$, satisfying $|\lambda_k|\le 1$ and whenever $|\lambda_k|=1$, then $\lambda_k$ is a simple eigenvalue. A more general condition than the latter one can be stated as follows: the zero solution of Eq. (\ref{1}) is stable if and only if the matrix $A$ has all eigenvalues $\lambda_k$, $k=1,2,\dots,n$, satisfying $|\lambda_k|\le 1$ and whenever $|\lambda_k|=1$, then all Jordan blocks corresponding to $\lambda_k$ is of size $1\times 1$. 

\begin{example}
Suppose $A$ is a $d\times d$-matrix such that $\sigma (A) \cap \Gamma = \{1\}$. Then, the matrix $A$ determines an operator from $\C^d\to \C^d$ that is polynomially bounded. In fact, if $\nu+1$ is the largest size of the Jordan blocks corresponding to the eigenvalue $1$ in the Jordan normal form of the matrix $A$, then we can split $\C^d$ as a direct sum $\C^d=E_1 \oplus E_2$, where $E_1$, $E_2$ are invariant under the linear transformation determined by $A$ such that $\sigma (A|_{E_1} )= \{1\}$ and $\sigma (A|_{E_2}) \subset \{ z\in \C : |z| < 1\}$. Using this splitting we can easily show that there exists a positive constant $M$ such that
\begin{align*}
   \| A^n\| \le Mn^\nu, \quad n\in \N.
\end{align*}
By Theorem \ref{the kat ext},
\begin{align*}
   \lim_{n\to \infty} \frac{1}{n^\nu} (A-I)^{\nu+1}A^n =0.
\end{align*}
This means that for all $x_0\in Range ((A-I)^{\nu+1})$,
\begin{align*}
   \lim_{n\to \infty} \frac{1}{n^\nu} A^nx_0 =0,
\end{align*}
and the limit is uniformly in $x_0\in Range ((A-I)^{\nu+1})$.
\end{example}

We note that this result could be proved by an elementary method as the space $\C^d$ can be split into direct sum of invariant subspaces corresponding to the eigenvalues of $A$. We also note that this simple technique is no longer valid in the infinite dimensional case when $A$ is a bounded linear operator in a (complex) Banach space $X$ since no splitting of $\X$ is guaranteed, so more advanced (abstract) techniques are needed. The reader is referred to \cite{arebathieneu} for more details about how to treat this case.

\begin{example}
Suppose that $\X$ is an infinite dimensional Banach space, $T$ is a contraction in $\X$ such that $\sigma (T)\cap \Gamma \subset \{1\}$,
and $K$ is a nilpotent bounded operator on $\X$ (that is, there exists $\nu\in\N$ such that $K^\nu=0$) that is commutative with $T$. Then, consider $B:=T+K$ and the asymptotic behavior of solution to the equation
\begin{align}\label{2}
   x(n+1)=Bx(n), \quad n\in \N.
\end{align}
Clearly, by \cite[Theorem 11.23 p. 293]{rud} the operator $B$ satisfies $\sigma(B)=\sigma (T+K)\subset \sigma (T)+ \sigma (K) =\sigma (T)+ \{0\}$, $\sigma (B)\cap \Gamma \subset \{1\}$. Moreover,
\begin{align*}
   \| B^n\|
    = \|(T+K)^n\| 
   &= \biggl\| \sum_{j=0}^{\nu -1} \frac{n(n-1)\cdots (n-j+1)}{j!}K^j T^{n-j}\biggr\| \\
   &\le \sum_{j=0}^{\nu -1} \frac{n(n-1)\cdots (n-j+1)}{j!}\| K\| ^j \\
   &\le Mn^{\nu-1},
\end{align*}
where $M$ is a positive number independent of $n$. Therefore, applying Theorem \ref{the kat ext} gives
\begin{align*}
   \lim_{n\to\infty} \frac{1}{n^{\nu-1}} (B-I)^{\nu} B^n=0 .
\end{align*}
In other words, $\lim_{n\to \infty} (1/n^{\nu-1}) B^nx_0=0$ uniformly in $x_0\in Range ((B-I)^{\nu})$.
\end{example}

When $\nu=1$, that is $K=0$, in the famous work \cite{kattza} by using Harmonic Analysis the authors proved that $\lim_{n\to \infty} (T-I)T^n=0$. An elegant proof of this result is given in \cite{vu1} in which a quotient space is employed to reduce the problem to an application of a Gel'dfand's theorem on spectra of isometries in Banach spaces. An elementary method of proof is offered in \cite{min} that is a basis for this paper for further extensions.

\begin{example}
Consider fractional difference equations of the form
\begin{align}\label{frac eq3}
   \Delta^\alpha x(n) =Ax(n)+y(n), \quad n\in \N,
\end{align}
where $x(n)\in \C^d$, $\alpha =1/2$, $A$ is a nilpotent $d\times d$-matrix such that $A^2 =0$
and $\{y(n)\}_{n=1}^\infty \in c_0(\C^d)$. Then, the condition $\Sigma \subset \{1\}$ means that 
\begin{align}\label{spec cond}
   \det (z^{1/2}(z-1)^{1/2} I-A) 
   \neq 0 \quad  \mbox{for all} \ z\in \Gamma \backslash \{1\}.
\end{align}
Since $\sigma (A)=\{0\}$, this condition is always satisfied. We check the boundedness of the resolvent $\{S_\alpha (n)\}_{n=1}^\infty$. By Theorem \ref{thm 4.9} (ii) we have
\begin{align*}
   S_\alpha (n) 
   =\sum_{j=0}^n \frac{\Gamma(n-j+(j+1)\alpha )}{\Gamma (n-j+1)\Gamma (j\alpha +\alpha)} A^j ,
\end{align*}
and hence, for each $n\in \N$, we have
\begin{align*}
   \| S_{1/2} (n) \| 
   &\le \frac{\Gamma(n+1/2 )}{\Gamma (n+1)\Gamma (1/2)}  + \frac{\Gamma(n )}{\Gamma (n)\Gamma (1)} \|A\|  \\
   &= \frac{(2n-1)!!\sqrt{\pi}/2^n}{n!\sqrt{\pi}}  + \|A\| \\
   &= \frac{(2n-1)!!}{(2n)!!}  + \|A\| \\
   &\le 1+\|A\| .
\end{align*}
By Theorem \ref{the kat for frac},
\begin{align}\label{final}
   \lim_{n\to \infty} [S_{1/2} (n+1) -S_{1/2}(n)]=0.
\end{align}
\end{example}

\begin{remark}
This result can be extended to the more general case when all eigenvalues $\lambda_k$ of $A$ has modulus less than $1$ such that condition (\ref{spec cond}) holds. In fact, in this case there exist constants $C>0$ and $0\le q<1$ such that $\| A^j\| \le Cq^j$, $j=1,2,\dots,$ so
\begin{align*}
   \| S_{1/2} (n) \| 
   \le \sum_{j=0}^n Cq^j 
   \le \frac{C}{1-q}, \quad  n\in \N .
\end{align*}
And (\ref{final}) again follows from Theorem \ref{the kat for frac}.
\end{remark}

\bibliographystyle{amsplain}

\end{document}